**Remarks on smoothness and finite-time blowup
solutions of the incompressible Navier-Stokes equation**

KAMAL N. SOLTANOV

Abstract. This article examines the smoothness of the solution to the Navier-Stokes equation from a novel perspective. Here, the existence of the smoother solution relative to x and to the time t was shown only for a finite time. Moreover, for each considered case of the problem, a blowup time for its solutions can be demonstrated.

1. Introduction

This article examines the well-known issue of the smoothness of solutions to the incompressible Navier-Stokes equation in the three-dimensional case. To study the smoothness of solutions to the PDE, various approaches were used, depending on the known functional inequalities related to the embedding theorems between functional spaces. Unfortunately, these approaches were unable to provide a solution to this problem. Therefore, it was necessary to find other approaches to this problem. We would note that this article uses a different approach.

In this article, we study the properties of solutions to the incompressible Navier-Stokes equation. Namely, properties of solutions of the following Cauchy problem for the well-known incompressible Navier-Stokes equation are investigated

1(1.1) $\qquad \frac{du}{dt} + \sum_{j=1}^{3} u_j \frac{\partial u}{\partial x_j} = \nu \Delta u - \nabla p + f, t > 0, \ x \epsilon R^3,$

(1.2) $\qquad divu = \sum_{i=1}^{3} \frac{\partial u}{\partial x_i} = 0, \ t > 0$

with initial conditions

(1.3) $\qquad u(0,x) = u_0(x), \ x \in R^3.$

Here, $u_0(x)$ is a given, $C^\infty$ divergence–free vector field on $R^3$, $f_i(t,x)$ are the components of a given, externally applied force (e.g., gravity), $\nu$ is a positive coefficient (the viscosity), and $\Delta u = \sum_{i=1}^{3} \frac{\partial^2}{\partial x_i^2}$ is the Laplacian in the space variables.

Well-known that the existence of weak solutions to this problem was proved by Leray and Hopf (see [13, 9, 10], and also e.g., [16, 30, 4, 7, 8, 29], which defined the space of weak solutions corresponding to the data spaces for this problem.

---

Presently, it has also been proven that the uniqueness of the weak solution to this problem holds under a few smoothness conditions on the data spaces of this problem (see [20]). We should note that these data spaces are everywhere dense in the appropriate spaces used under the proof of the existence theorem (see [20]). It is clear that due to the uniqueness of the weak solution, under conditions of the uniqueness theorem, all possible solutions to this problem are unique.

We should note that the regularity and partial regularity of the weak solutions to the Navier-Stokes equation have been studied and are also studied now by a few authors in articles (see e.g., [1, 2, 3, 4, 6, 11, 12, 13, 14, 15, 16, 17, 18, 19, 28, 30, 31, 32, 33], their references), in which the mentioned smoothness were studied under different complementary conditions in various ways.

This article examines the smoothness of the solution to the Navier-Stokes equation from a different perspective.

This article is organized as follows. Section 2 provides the main results, certain known solvability theorems, and proves necessary auxiliary results. In Section 3, the existence of a finite time of a strong solution to the considered problem is established under some sufficient conditions. In Section 4, the existence of a finite time of a smoother solution to the problem relative to $x$ is established under some sufficient conditions. finite time. In Section 6, the result of this type, as the above result for $t$ and also for $x$, was shown.

## 2. Preliminary & Main Results

As is known, to study the smoothness of solutions to problems in PDE theory, the known embedding theorems between functional spaces are usually applied. We will also use these theorems in another way. However, we will use a different approach to study this problem. In this case, we will use the embedding of the space $H^2(R^3)$ intothe space $L^\infty(R^3)$, but elsewhere.

We assume the data functions sufficiently smoothly, as we will study the existence of a smooth solution to the Navier-Stokes equation. We introduce the following notations $Q = R_+ \times R^3, H = \left(L^2(R^3)\right)^3, H^m = W^{m,2}(R^3), L^2\left(R_+; \left(L^2(R^3)\right)^3\right) = L^2(R_+; H)$.

Let $u_0 \in (H^m) \cap V(R^3)$ and $f \in W^{k,2}(R_+; H^m))$, where $m \geq 1$ and $V = V(R^3)$ is the space of the Leray weak solutions. Let $A(D_x)$ be linear differential operator acting in $H$ such that $KerA = \{0\}$, $Im(A(D_x))$ everywhere dense in $H$, which has the form
$A(D_x) = \sum_{|\alpha| \geq 0}^m a_\alpha, D_x^\alpha = D_{x_1}^{\alpha_1} \cdots D_{x_n}^{\alpha_n}, \alpha = (\alpha_1, \cdots, \alpha_n), |\alpha| = \sum_{i=1}^n \alpha_i, \ \alpha_i \geq 0$,
where $a_\alpha \geq 0, \ 0 \leq |\alpha| \leq m$ are constants, therefore, $A(D_x)$ commutes with $D_t, D_{x_i}, \Delta, \nabla$ and $\alpha_i, k, m \in N$, and $D_t = \frac{\partial}{\partial t}, D_{x_i} = \frac{\partial}{\partial x_i}$. Moreover, assume $A(D_x)$ and $A^{1/2}(D_x)$ is the self-adjoint operator, $D_t^k A(D_x)$ acting in $L^2(H)$.

We wish to note that for the study of the posed problem, we will consider this problem in the following formulation: $u$ and $p$ satisfy the equation

$$\langle \frac{\partial u}{\partial t} + (u \cdot \nabla)u - \mu \Delta u + \nabla p, v \rangle = \langle f, v \rangle,$$

for any $v$, at least, from $V(R^3)$ and a.e. $t \in (0,T)$, and fulfills (1.2), (1.3), where $T \in R_+$ will later be defined. Here $\langle \cdot, \cdot \rangle \equiv \int \cdot \times \cdot \, dx$.

We will now present the main results of this article. We assume that the dates ($u_0$ and $f$) possess the required smoothness, which is necessary for the considered case. In other words, the selection of the smoothness of the dates depends on the unknown order of smoothness of solutions, which we would like to show in the present case. As the formulations of the four theorems are similar, it will be best to present these together, since we have studied only four cases here.

**Theorem 1.** *Let $u \in L^2(V)$ is a weak solution to Problem 1, if $u_0 \in H^3, f \in L^2(R_+; H^2)$ then there exists a finite blowup time $T_1 > 0$ such that for each $T < T_1$ the solution $u(t,x)$ to Problem 1 belongs to the space*

$$L^\infty(0,T;H^1) \cap L^2(0,T;H^2 \cap V) \cap W^{1,2}(H),$$

*in other words, $u(t,x)$ is a strong solution $\left(consequently, p \in L^2(H)\right)$ and, $T_1 = T_1(\mu, u_0, f) > 0$ is defined with equality* (3.8).

**Theorem 2.** *Let $u \in L^2(V)$ is a weak solution to Problem 1 if $u_0 \in H^4, f \in L^2(R_+; H^4)$, then there exists a finite blowup time $T_2 > 0$ such that for each $T < T_2$ solution $u$ to the Problem 1 belongs to the space*

$$L^\infty(0,T;H^2) \cap L^2(0,T;H^3 \cap V) \cap W^{1,2}(0,T;H^1)$$

*in other words, $u$ is a smoother solution $\left(consequently, p \in L^2(H^1)\right)$, where $T_2 = T_2(\mu, u_0, f) > 0$ is defined with the equality of type (3.8)(since time $T_2$ is also found using Proposition 1 and the corollary of type Corollary 1).*

**Theorem 3.** *Let $u \in L^2(V)$ is a weak solution to Problem 1, if $u_0 \in H^4, f \in W^{1,2}(R_+; H^3)$, then there exists a finite blowup time $T_3 > 0$ such that for each $T < T_3$ the solution $u$ to the Problem 1 belongs to the space*

$$L^\infty(0,T;H^2) \cap L^2(0,T;H^2) \cap W^{1,2}(0,T;H^1), T < T_3$$

*In other words, $u$ is a smoother solution concerning $t$ (consequently, $p \in L^2(H^1)$), where $T_3(\mu, u_0, f) > 0$ satisfies the equality of type (3.8)(since time $T_3 > 0$ is also found using Proposition 1 and the corollary, of type Corollary 1).*

**Theorem 4.** *Let $u \in L^2(V)$ is a weak solution to Problem 1, if $u_0, u_1 \in H^4, f \in W^{2,2}(R_+; H^4)$, then there exists a finite blowup time $T_4(\mu, u_0, u_1, f) > 0$, such that for each $T < T_4$ the solution $u$ to Problem 1 belongs to the space*

$$W^{1,\infty}(0,T;H^2) \cap W^{2,2}(0,T;H \cap V) \cap W^{1,2}(0,T;H^2), T < T_4$$

*in other words, $u$ is a smoother solution concerning $t$ (consequently, $p \in L^2(H^1)$), where $T_4 > 0$ satisfies the equality of type (3.8)(since time $T_4$ is found also using Proposition 1 and the corollary, of type Corollary 1).*

To continue these studies, it is necessary to present certain known results on solvability and to prove some auxiliary results. We begin with the special case of two general solvability theorems, which were proved in articles [22, 23, 24, 25, 26, 27] (see also the book [21]); therefore, we provide them without their proofs.

Let $X, Y$ be a reflexive Banach spaces and $X^*, Y^*$ *are* their dual spaces, let $Y$ be reflexive with strictly convex norm, together with $Y^*$ (well-known, this condition is not a complementary condition see [5]), and $f: D(f) \subseteq X \rightarrow Y$ be an operator.

Study of the following equation

$$f(x) = y, x \in D(f), y \in Y$$

Consider the following conditions. Let the closed ball $B_{r_0}^X(0)$ of $X$ belongs to $D(f) \subseteq X$ and on ball $B_{r_0}^X(0)$ the following conditions are fulfilled

(1) $f: B_{r_0}^X(0) \longrightarrow Y$ be the continuous operator that is bounded on $B_{r_0}^X(0)$, i.e. $\|f(x)\|_Y \leq \mu(x_X)$, $\forall x \in B_{r_0}^X$, $\mu: R_+ \to R_+$ $\mu \in C^0$;

(2) There is a mapping $g: D(g) \subseteq X \longrightarrow Y^*$ that $D(f) \subseteq D(g)$ and for any $S_r^X(0) \subset B_{r_0}^X(0), 0 < r \leq r_0, cl_{Y^*} g(S_r^X(0)) \equiv \overline{g(S_r^X(0))} = S_r^{Y^*}(0), S_r^X(0) \subseteq g^{-1}\left(S_r^{Y^*}(0)\right)$

$$\langle F(x), g(x)\rangle \geq \rho(\|x\|_X), \quad a.e.\ x \in B_{r_0}^X(0), \quad \rho(r_0) \geq \delta_0 > 0$$

holds[1], $\frac{\rho(\tau)}{\tau} \nearrow \infty$ for $\tau \nearrow \infty$, where $\rho: R_+ \to R$, $\rho, \mu$ are continuous functions; moreover, the function $\lambda$ is such that $\rho(\tau) \geq 0, for\ \tau \geq \tau_0 \geq 0$, and increases for $r_0 > \tau \geq \tau_0, \delta_0, \tau_0 \geq 0$ are constants;

(3) for almost every $\tilde{x} \in Int B_{r_0}^X(0)$ possesses a neighborhood $U_\varepsilon(\tilde{x}), \varepsilon \geq \varepsilon_0 > 0$ that the following inequality

$$\| f(x_1) - f(x_2) \|_Y \geq \varphi(\| x_1 - x_2 \|_X, \varepsilon, x), x_1, x_2 \in U_\varepsilon(\tilde{x}) \cap B_{r_0}^X(0),$$

holds, where $\varphi(\tau, \tilde{x}, \varepsilon) \geq 0$ is a continuous function at $\tau$ and $\varphi(\tau, \tilde{x}, \varepsilon = 0 \leftrightarrow \tau = 0)$ (in particular, $\varepsilon = \varepsilon_0 = r_0, \tilde{x} = 0 \to U_\varepsilon(0) = B_{r_0}^X(0)$, consequently, $\varphi(\tau, \tilde{x}, \varepsilon) = \varphi(\tau, 0, r_0)$ on $B_{r_0}^X(0)$.

**Theorem 5.** *Let $X, Y$ be Banach spaces such as above, and $f : D(f) \subseteq X \to Y$ is a continuous operator. Assume that on the closed ball $B_{r_0}^X(0) \subseteq X$ conditions (1) and (2) are fulfilled, then the image $f\left(B_{r_0}^X(0)\right)$ of the ball $B_{r_0}^X(0)$ belongs to an everywhere dense subset of the bodily subset $M \subseteq Y$, that has the form $M = \{y \in Y | \langle y, g(x)\rangle \leq \langle f(x), g(x)\rangle, \forall x \in S_{r_0}^X(0) \}$*

*Furthermore, if in addition $f\left(B_{r_0}^X(0)\right)$ closed subset of $Y$ or the condition (3) is fulfilled, then the image $f\left(B_{r_0}^X(0)\right)$ of the ball $B_{r_0}^X(0)$ is a bodily subset of $Y$. Moreover, $f\left(B_{r_0}^X(0)\right)$ contains the above bodily subset $M \subseteq Y$* [2].

We should note that, under the application of Theorem 5, if it is necessary, one can begin to use the elliptic regularization.

Let $X$ and $Y$ $\overline{be}$ Banach spaces with duals $X^*$ and $Y^*$, respectively, $Y$ be a reflexive Banach space, $M_0 \subseteq X$ be a weakly complete "reflexive" $pn$−space (see e.g. [21, 26]), $X_0 \subseteq M_0 \cap Y$ be a separable vector topological space such that $cl_{M_0} X_0 \equiv M_0$, $cl_Y X_0 \equiv Y$ Consider the following problem:

---

[1] In particular; a mapping $g \equiv L$ is a bounded linear operator as $g \equiv L: X \longrightarrow Y^*$ that satisfies the conditions of *(ii)*.

[2] We should note that the works mentioned above contain more general sufficient conditions for the closure of the image $f\left(B_{r_0}^X(0)\right)$.

In particular, we can assume the condition that operator $f$ possess property $P$, i.e., $f: D(f) \subseteq X \to Y$ iff each precompact subset $U \subseteq Y$ from $Im(f)$ has a such (general) subsequence $U_0 \subseteq U$ that exists a precompact subset $G \subseteq X$ that satisfies the inclusion $f^{-1}(U_0) \subseteq G$ and $U_0 \subseteq f(G \cap D(f))$

$$\frac{dx}{dt} + f(t, x(t)) = y(t), y \in L^{p_1}(0, T; Y), x(0) = 0 \quad (2.1)$$

Let the following conditions be fulfilled: *i)* $f: \mathbb{P}_{1p_0p_1}(0, T; M_0, Y) \to L^{p_1}(0, T; Y)$ is a weakly compact mapping[3], where

$$\mathbb{P}_{1p_0p_1}(0, T; M_0, Y) \equiv L^{p_0}(0, T; M_0) \cap W^{1,p_1}(0, T; Y) \cap \{x(t)|x(0) = 0\}$$

$$1 < max\{p_1, p_1'\} \leq p_0 < \infty, p_1' = \frac{p_1}{p_1 - 1}$$

*(ii)* There is a linear, continuous operator $L: W^{s,p_2}(0, T; X_0) \to W^{s,p_2}(0, T; Y^*), s \geq 0, p_2 \geq 1$ such that $L$ commutes with $\frac{d}{dt}$ and the conjugate operator $L^*$ has $ker(L^*) = \{0\}$

*(iii)* There exists a continuous function $\varphi: R_+ \cup \{0\} \to R$ and numbers $\tau_0 \geq 0$ and $\tau_1 > 0$ such that $\varphi(r)$ is nondecreasing for $\tau_1 \geq \tau_0, \varphi(\tau_1) > 0$ and operators $f$ and $L$ satisfy the following inequality for any $x \in L^{p_0}(0, T; X_0)$

$$\int_0^T \langle f(t, x(t)), Lx(t) \rangle dt \geq \varphi([x]_{L^{p_0}(M_0)})[x]_{L^{p_0}(M_0)};$$

*(iv)* There exists a linear bounded operator $L_0: X_0 \to Y$ and constants $C_0 > 0, C_1, C_2 \geq 0, \theta > 1$ such that the inequalities

$$\int_0^T \langle \xi(t), L\xi(t) \rangle dt \geq C_0 \|L_0 \xi\|_{L^{p_1}(0,T;Y)}^\theta - C_1,$$

$$\int_0^t \langle \frac{dx}{d\tau}, Lx(\tau) \rangle d\tau \geq C_0 \|L_0 x\|_Y^\theta(t) - C_2, \quad a.e. t \in ]0, T]$$

hold for any $x \in W^{1,p_0}(0, T; X_0)$ and $\xi \in L^{p_0}(0, T; X_0)$.

**Theorem 6.** *Assume that conditions (i) - (iv) are fulfilled. Then the Cauchy problem (2.1) is solvable in $\mathbb{P}_{1p_0p_1}(0, T; M_0, Y)$ in the following sense*

$$\int_0^T \langle \frac{dx}{dt} + f(t, x(t)), y^*(t) \rangle dt = \int_0^T \langle y(t), y^*(t) \rangle dt, \quad \forall y^* \in L^{p_1'}(0, T; Y^*),$$

*for any $y \in G \subseteq L^{p_1}(0, T; Y)$, where $G = \cup_{r \geq \tau_1} G_r$ that is defined as*
$G_r \equiv \left\{ (y \in L^{p_1}(0, T; Y) \Big| \int_0^T |\langle y, Lx \rangle| dt \leq \int_0^T \langle f(t, x), Lx \rangle dt - c \right) \forall x \in L^{p_0}(0, T; X_0),$
$[x]_{L^{p_0}(M_0)} = r \right\}, \quad C_2 < c < \infty.$

Now we will study the following Cauchy problem, which is the auxiliary result.

**Proposition 1**. *Let the function $y(t)$ is the solution to the following Cauchy problem*

$$\frac{dy}{dt} \leq ay^2 + b, \quad y(0) = y_0$$

*then $y(t)$ satisfies the following inequality*

$$y(t) \leq \frac{y_0 + \left(\frac{a}{b}\right)^{-\frac{1}{2}} \tan\left((ab)^{\frac{1}{2}}t\right)}{1 - \left(\frac{a}{b}\right)^{\frac{1}{2}} y_0 \tan\left((ab)^{\frac{1}{2}}t\right)},$$

---

[3] See [26]. This condition is explained in a similar condition 1 of Appendix B.

where $a, b > 0$ are constants.

Proof. The proof follows from the studied Cauchy problem for the Riccati equation. Indeed, we have

$$\frac{dy}{dt} \leq ay^2 + b \rightarrow \left(\frac{a}{b}\right)^{1/2} \frac{dy}{dt} \leq (ab)^{1/2} \left(\frac{a}{b}y^2 + 1\right) \rightarrow$$

Let $z = \left(\frac{a}{b}\right)^{1/2} y$ then we get

$$\frac{dz}{dt} \leq (ab)^{1/2}(z^2 + 1) \rightarrow \tan^{-1} z(t) \leq \tan^{-1} z(0) + (ab)^{1/2}$$

$$\rightarrow \tan^{-1} \left(\frac{a}{b}\right)^{1/2} y(t) \leq \tan^{-1} \left(\frac{a}{b}\right)^{1/2} y(0) + (ab)^{1/2} \rightarrow$$

$$y(t) \leq \left(\frac{a}{b}\right)^{-1/2} \tan\left(\tan^{-1} \left(\frac{a}{b}\right)^{1/2} y(0) + (ab)^{1/2} t\right) =$$

$$= \left(\frac{a}{b}\right)^{-1/2} \frac{\left(\frac{a}{b}\right)^{1/2} y_0 + \tan\left((ab)^{1/2} t\right)}{1 - \left(\frac{a}{b}\right)^{1/2} y_0 \tan\left((ab)^{1/2} t\right)} \rightarrow$$

$$\rightarrow y(t) \leq \frac{y_0 + \left(\frac{a}{b}\right)^{-1/2} \tan\left((ab)^{1/2} t\right)}{1 - \left(\frac{a}{b}\right)^{1/2} y_0 \tan\left((ab)^{1/2} t\right)}.$$

Consequently, the obtained estimation has a sense only for a finite time, i.e., it must fulfill the following inequality

$$t < (ab)^{-1/2} \tan^{-1}\left(\left(\frac{a}{b}\right)^{-1/2} y_0^{-1}\right).$$

## 3. On the existence of the strong solution

Here, we will introduce the main estimates of this article.

First, we assume that $A(Dx)u = (-\Delta + I)u$, and consider the following equation

$$\left\langle \frac{\partial u}{\partial t} + \sum_{j=1}^{3} u_j \frac{\partial u}{\partial x_j} - \nu \Delta u + \nabla p - f, -\Delta u + u \right\rangle = 0$$

hence follows

$$\frac{1}{2}\frac{d}{dt}(\|u\|^2 + \|\nabla u\|^2) + \nu\|\Delta u\|^2 + \nu\|\nabla u\|^2 + \langle((u \cdot \nabla)u), -\Delta u + u\rangle$$
$$= \langle f, -\Delta u + u\rangle, \quad (3.1)$$

due to the conditions, the problem (see also (1.2)), we have $-\langle \nabla p, \Delta u\rangle$

$$\langle \nabla p, u\rangle = \langle p, \Delta \nabla u\rangle - \langle p, \nabla u\rangle = 0.$$

As we will study the smoothness of the solution, here and in what follows, we will assume that the function $u \in V(R^3)$ is sufficiently smooth, more exactly, $u \in C^0(Q) \cap W^{m,p}(Q)$, since this class of functions is everywhere dense, e.g., in spaces $L^p(Q), p \geq 1$. Moreover, the norm of the space $(L^2(Q)) \equiv H^3$, we will write as
$\|\cdot\|$ for brevity.

3.0.1. *Necessary Estimations.* As is known, the weak solution to the Navier-Stokes equation belongs to the space $L^\infty(R_+; H) \cap L^2(R_+; V)$.

To obtain the a priori estimates and details of the receipt of the necessary inequalities, we will rewrite equality 1 in the following form:

$$\int \left( \frac{\partial u}{\partial t} + \sum_{j=1}^{3} u_j \frac{\partial u}{\partial x_j} - \nu \Delta u + \nabla p \right)(-\Delta u + u) dx = \int f(-\Delta u + u)\, dx$$

hence we get

$$\int \frac{1}{2} \frac{\partial}{\partial t}(|u|^2 + |\nabla u|^2) dx + \nu \|\Delta u\|^2 + \nu \|\nabla u\|^2 = \int [((u \cdot \nabla)u) + \nabla p - f](\Delta u - u)\, dx$$

$$\leq \int ((u \cdot \nabla)u)\Delta u dx + \frac{C}{2\nu}\|f\|_{L^2(H^2)}^2 + \frac{\nu}{2C}\|\Delta u\|^2 + \frac{\nu}{2}\|\nabla u\|^2 \leq \int \frac{c_1}{2\nu}(|u|^2 + |\nabla u|^2)^2 dx$$

$$+ \frac{\nu}{2c_1}\|\Delta u\|^2 + \frac{C}{2\nu}\|f\|_{L^2(H^2)}^2 + \frac{\nu}{C_1}\|\Delta u\|^2 + \frac{\nu}{2c_2}\|\nabla u\|^2$$

according to the smoothness of the dates of the problem, where $C, C_1, C_2, c_1, c_2 \geq 1$, $\langle \nabla p, (\Delta u - u) \rangle = 0$, and also $\langle (u \cdot \nabla)u, u \rangle = 0$. Thus, we arrive at the inequality

$$\frac{\nu}{2}(\|\Delta u\|^2 + \|\nabla u\|^2) + \int \frac{c}{\nu} \frac{\partial}{\partial t}(|u|^2 + |\nabla u|^2) dx \leq \int \frac{1}{8\nu}(|u|^2 + |\nabla u|^2)^2 dx + \frac{C}{2\nu}. \quad (3.2)$$

Whence, follows that it is necessary to study the following Cauchy problem for the ordinary differential equation, as the first two terms on the left side are positive

$$\frac{1}{2}\frac{d}{dt}(|u|^2 + |\nabla u|^2) \leq \frac{1}{8\nu}(|u|^2 + \nabla u^2)^2 + \frac{C}{\nu^2}\|f\|^2, \ u(0,x) = u_0, \nabla u(0,x) = \nabla u_0 \quad (3.3)$$

*To study this problem, we will use Proposition* 1; *from here, the result follows.*

[4]Here, one can use, $\|f\|_{L^2(R_+)}$, but for brevity, we use the norm of $f$.

**Corollary 1.** *The expression satisfies* $|u|^2 + |\nabla u|^2$ *the following estimation*

$$|u|^2 + |\nabla u|^2 \leq$$

$$\leq \frac{|u_0|^2 + \nabla u_0{}^2 + \left(\frac{2}{C}\|f\|^2\right)^{1/2} \tan\left(\left(\frac{C}{8\nu^2}\|f\|^2\right)^{1/2} t\right)}{1 - (|u_0|^2 + \nabla u_0{}^2)\left(\frac{2}{C}\|f\|^2\right)^{-1/2} \tan\left(\left(\frac{C}{8\nu^2}\|f\|^2\right)^{1/2} t\right)} = M \quad (3.4)$$

where $M = M(t, u_0(x), f, \nu)$. Where does it imply that this inequality makes sense only if the following inequality

$$\left(\frac{C}{8\nu^2}\|f\|^2\right)^{-1/2} \tan^{-1}\left((|u_0|^2 + \nabla u_0{}^2)\left(\frac{2}{C}\|f\|^2\right)^{1/2}\right)^{-1} < t \quad (3.5)$$

holds.

So, we can continue to the receipt of the necessary a priori estimate for the problem. Consider the inequality

$$\frac{1}{2}\frac{d}{dt}(|u|^2 + |\nabla u|^2) + \frac{\nu}{2}|\Delta u|^2 + \frac{\nu}{C}|\nabla u|^2 \leq$$

$$\leq \frac{1}{2}(|u|^2 + |\nabla u|^2)^2 + \frac{C}{2\nu}\|f\|_{L^2(H^2)}^2$$

Then, due to Corollary 1, we get

$$\frac{d}{dt}(\|u\|^2 + \|\nabla u\|^2) + \nu\|\Delta u\|^2 + \nu c_1\|\nabla u\|^2 \leq$$
$$\frac{1}{\nu}\int M(t,u_0(x),\cdots)^2 dx + C\left(\|f\|_{L^2(H^2)},\nu\right). \quad (3.6),$$

and if we integrate this inequality concerning $t$ we obtain

$$(\|u\|^2 + \|\nabla u\|^2)(t) + \int_0^t \left(\frac{\nu}{2}\|\Delta u\|^2 + \nu c_1\|\nabla u\|^2\right) ds \leq$$
$$\leq \int_0^t \int \frac{1}{\nu} M(s,u_{0,\cdots})^2 dx ds + C\left(\|f\|_{L^2(H^2)},\nu,\|u_0\|_{H^1}\right) \quad (3.7)$$

From here follows, if $u_0 \in H^3$ and $f \in L^2(H^2)$, then the right side of the inequality (3.7) is bounded, and consequently, the left side will also be bounded.

So, we get that $u \in L^\infty(0,T;H^1) \cap L^2(0,T;H^2)$, and consequently,
$$B(u,u) = (u \cdot \nabla)u \in L^2(0,T;H)$$

For each $T \in (0,T_1)$, where $T_1$ fulfills the equation

$$T_1 = \left(\frac{C}{4\nu^2}\|f\|^2\right)^{-1/2} \tan^{-1}\left((|u_0|^2 + |\nabla u_0|^2)\left(\frac{2\|f\|^2}{C}\right)^{-1/2}\right)^{-1}. \quad (3.8)$$

To prove the existence theorem for the studied problem, we will use Theorem 5. Thus, if we denote by $F$ the operator generated by the studied problem and the operator $g(v) = A(Dx)v = (-\Delta + I)v$, then, since the boundedness of the operator $F$ follows from the above result, we must show that

$$\langle F(u), g(u) \rangle \geq \rho(\|u\|)$$

holds for any $u \in L^\infty(0,T;H^1) \cap L^2(0.T;H^2) \cap L^2(0,T;L^\infty)$, where $\lambda: R_+ \to R, T < T_1$, and satisfies the condition of Theorem 5.

To prove the necessary inequality, one can use the inequalities mentioned above. Therefore, consider the following

$$\langle F(u), g(u) \rangle = \langle \left[\frac{\partial u}{\partial t} + ((u \cdot \nabla)u) - \nu\Delta u + \nabla p\right], (-\Delta u + u)\rangle$$
$$= \frac{1}{2}\frac{d}{dt}(\|u\|^2 + \|\nabla u\|^2) + \nu(\|\Delta u\|^2 + \|\nabla u\|^2) + \langle ((u \cdot \nabla)u), (-\Delta u)\rangle$$
$$+ \langle \nabla p, (-\Delta u + u)\rangle =$$

taking into account the condition of the problem, we get

$$= \frac{1}{2}\frac{d}{dt}\int (|u|^2 + |\nabla u|^2)dx + \nu(\|\Delta u\|^2 + \|\nabla u\|^2) + \int ((u \cdot \nabla)u)(-\Delta u)dx$$

Hence follows

$$\langle F(u), g(u)\rangle \geq \frac{1}{2}\frac{d}{dt}(\|u\|^2 + \|\nabla u\|^2) + \frac{\nu}{2}\|\Delta u\|^2 + \nu\|\nabla u\|^2$$
$$- \frac{1}{2\nu}\int M(\|f\|,u_0,\cdots)^2 dx$$

from here by integrating at $t$, we get

$$\int_0^t \langle F(u), g(u)\rangle ds \geq \frac{1}{2}(\|u\|^2 + \|\nabla u\|^2)(t) + \int_0^t \left(\frac{\nu}{2}\|\Delta u\|^2 + \nu\|\nabla u\|^2\right)(s)\,ds$$
$$- \frac{1}{2\nu}\int_0^t \left(M(\|f\|,u_0,\nu,\cdots)\right)^2 ds$$

Thus, we arrive
$$\int_0^T \langle F(u), g(u)\rangle(t)\, dt \geq \hat{c}\|u\|^2_{L^2(H^2)\cap L^\infty(H^2)} - \int_0^T |M_1(u_o, f)|^2 dt$$

according to Corollary 1.

We should note that in this case, the image $F(U)$, of each convex subset $U \subset X$ is a convex set; therefore, to show closure of this image, the Mazur Theorem is used. The a priori estimates and the above inequality show that under given conditions on the initial dates $u_0$ and the externally applied force $f$, one can prove the existence of the smooth solution using Theorem 5. We should note that one can apply to the considered problem of Theorem 5 after using the elliptic regularization.

Moreover, it needs to be noted that the following equality is fulfilled due to the obtained result on the existence of a solution

$$\langle \frac{\partial u}{\partial t}, v\rangle = \langle f + \Delta u - ((u\cdot\nabla)u), v\rangle$$

for any $v \in L^2(H)$. Moreover, it satisfies the following inclusion $(u\cdot\nabla)u \in L^2(H)$. Really, for $u \in L^2(0, T; H^2) \cap L^\infty(0, T; H^1)$ we have

$$\int\int |(u\cdot\nabla)u|^2\, dxdt \leq c\int \|u\|^2 \|\nabla u\|^2 dt \leq c\|u\|^2_{L^2(L^\infty)} \|\nabla u\|^2_{L^\infty(H)},$$

according to the Sobolev inequality, consequently, $(u\cdot\nabla)u \in L^2(0, T; H))$. Then we get that $u_t \in L^2(0, T; H))$ is true, according to the above equation, since all other components belong to this space. In addition, as this equation must contain the derivative of the pressure, then $\nabla p \in L^2(0, T; H))$ is also valid.

Consequently, we have proved that under our conditions on the initial value and the outer force, the Navier-Stokes equation has a strong solution. Now we will show that the solution to Problem 1 can be smoother.

## 4. Other results on smoothness

Let the dates of the considered problem be sufficiently smooth. Consider the following expression to study the additional property of the solution to the studied problem with the $A(D_x) = u + \Delta^2 u$

$$0 = \langle \frac{\partial u}{\partial t} + (u\cdot\nabla)u - \nu\Delta u + \nabla p - f, u + \Delta^2\rangle =$$

$= \frac{1}{2}\frac{d}{dt}(\|u\|^2 + \|\Delta u\|^2) + \langle (u\cdot\nabla)u, u + \Delta u^2\rangle + +\nu\|\nabla\Delta u\|^2 + \nu\|\nabla u\|^2 + \langle \nabla p, u + \Delta^2 u\rangle - \langle f, u + \Delta^2 u\rangle$ (4.1)

To study the smoothness of the solution in this case, we will use the previous approach. Consequently, we must receive the necessary a priori estimates. So, using the same argument as above, it is enough to study terms containing $\nabla p$ and $(u\cdot\nabla)u$. We can note that $\langle \nabla p, u + \Delta^2 u\rangle = 0$ according to the conditions of the considered Problem. Therefore, remains to study the term $\langle (u\cdot\nabla)u, \Delta^2 u\rangle$ that needs to be estimated.

Since we will use the above approach, we will rewrite the equality (4.1) in a different form

$$\int \left[\left(\frac{1}{2}\frac{\partial}{\partial t}(|u|^2 + |\Delta u|^2)\right) + \nu(|\nabla u|^2 + |\nabla \Delta u|^2)\right] dx + \int ((u \cdot \nabla)u)\Delta^2 u\, dx = \int f(u + \Delta^2 u)dx \quad (4.2).$$

Now, we will investigate $\int ((u \cdot \nabla)u(\Delta^2 u))dx$, therefore, we consider it separately

$$|\langle ((u \cdot \nabla)u), \Delta^2 u \rangle| = |\langle \nabla((u \cdot \nabla)u), \nabla \Delta u \rangle|$$
$$\leq \int \frac{1}{2\nu}|\nabla((u \cdot \nabla)u)|^2 dx + \frac{\nu}{2}\|\nabla \Delta u\|^2$$

Consider the first term on the right-hand side of the above inequality

$$\int (\nabla(u \cdot \nabla)u)^2 dx = \int \left(((\nabla u \cdot \nabla)u + ((u \cdot \nabla)\nabla u))\right)^2 dx$$

Whence, it isn't difficult to see that the inequality

(4.3) $$\int (\nabla(u \cdot \nabla)u)^2 dx \leq \int (|u|^2 + |\Delta u|^2)^2 dx$$

holds. Takes into account the inequality (4.3) in the equality (4.2), and draws the estimates we have

(4.4) $$\int \frac{1}{2}\frac{\partial}{\partial t}(|u|^2 + |\Delta u|^2)dx + \frac{\nu}{2}\|\nabla \Delta u\|^2 + \frac{(2c-1)\nu}{2c}\|\nabla u\|^2 \leq$$
$$\leq \int \frac{2}{\nu}(|u|^2 + \Delta u^2)^2 dx + \frac{c}{2\nu}\|f\|_{H^4}^2$$

Thus, it is necessary to study the following Cauchy problem for the ordinary differential inequality

(4.5) $$\frac{1}{2}\frac{d}{dt}(|u|^2 + \Delta u^2) \leq \frac{2}{\nu}(|u|^2 + |\Delta u|^2)^2 + \frac{c}{2}\|f\|_{H^4}^2$$

obtained from the inequality (4.4) with the initial conditions $u(0) = u_0$ and $\Delta u(0) = \Delta u_0$ using Corollary 1 as in Section 3.

Not difficult to see that this inequality is of the same type as in the previous case, and it gives the estimation of the same type as mentioned above. Consequently, the obtained estimation is true for a finite time $T_2 > 0$ that can be defined in the same way as in Section 3.

Further, by continuing the argumentation, analogous to the arguments in the previous case, we will arrive at the necessary a priori estimates from which follows inclusion

$$u \in L^\infty(0, T; H^2) \cap L^2(0, T; H^3)$$

where follows $(u \cdot \nabla)u \in L^2(0, T; H^2), 0 < T < T_2$. We have

$$|D_x(u \cdot \nabla)u| \leq c(|D_x u||\nabla u| + |u||D_x \nabla u|) \leq c_1(|\nabla u||\nabla u| + |u||\Delta u|)$$

then, to show the necessary inclusion, we act as follows

$$\int_0^T \int (|\nabla u||\nabla u| + |u||\Delta u|)^2 dx\, dt \leq c \int_0^T (\|\nabla u\|^2\|\nabla u\|^2 + \|u\|^2\|\Delta u\|^2)dt \leq$$
$$\leq c\|\nabla u\|_{L^\infty(H)}^2 \nabla u_{L^2(L^\infty)}^2 + \|u\|_{L^\infty(L^\infty)}^2 \|\Delta u\|_{L^\infty(H)}^2$$

due to the Sobolev inequality, this shows that $(u \cdot \nabla)u \in L^2(0, T; H^1)$ is true, due to inclusion $u \in L^\infty(0, T; H^2) \cap L^2(0, T; H^3)$.

Now we can repeat the last part of the previous section.

The above investigations demonstrate that if we continue the study the same way as the previous cases, choosing the smooth dates for the problem, then we will obtain solutions that possess greater smoothness according to the data. Moreover, the used method shows that a hierarchy exists dependent on the choice of the operator $A(D_x)$.

## 5. Smoothness of solutions relative to the variable $t$

To study the properties of the solution relative to $t$ when the data relative to $t$ be smoother, we will use the operator $D_t A(D_x)$. So, we act in the following way

$$\langle \frac{\partial u}{\partial t} + (u \cdot \nabla)u - \nu\Delta u + \nabla p - f, u_t - \Delta u_t \rangle = 0$$

whence, as in the previous Section 3, we get

$$0 = \langle \frac{\partial u}{\partial t} + (u \cdot \nabla)u - \nu\Delta u + \nabla p - f, u_t - \Delta u_t \rangle \rightarrow$$

$$\int \left[ |u_t|^2 + |\nabla u_t|^2 + \frac{1}{2}\frac{d}{dt}(\nu|\Delta u|^2 + |u|^2) \right] dx$$

(5.1)  $\int [(u \cdot \nabla)u + \nabla p](u_t - \Delta u_t) dx = \int f(u_t - \Delta u_t) dx$

By the above reasons, due to the conditions we have $\int \nabla p(u_t - \Delta u_t) dx = 0$, consequently, it remains to estimate $\int (u \cdot \nabla)u(\Delta u_t) dx$. To estimate this, we act in the following way

$$\left| \int (u \cdot \nabla)u(\Delta u_t) dx \right| \leq \int |(\nabla(u \cdot \nabla)u)\nabla u_t| dx \leq$$

$$\frac{1}{2\nu} \int |[\nabla(u \cdot \nabla)u]|^2 dx + \frac{\nu}{2} \int |\nabla u_t|^2 dx$$

Consequently, it remains to estimate the integral $\int |(\nabla(u \cdot \nabla)u)|^2 dx$. It is not difficult to see that the following inequality is true

(5.2)  $$\int |(\nabla(u \cdot \nabla)u)|^2 dx \leq 4 \int (|\Delta u|^2 + |u|^2) dx,$$

according to the calculations conducted in Sections 3 and 4.

Thus, we obtain the Cauchy problem for the ordinary differential inequality, which is similar to the problem studied in Section 4. It is not difficult to see that the mentioned Cauchy problem is the same as in Section 4, i.e., the differential inequality is similar to (4.5), according to the expressions (5.1) and
(5.2).

Therefore, by repetition of the discussions analogous to the above sections, we can obtain the necessary a priori estimates for the problem considered here. To complete our study in this case, it's sufficient to use analogous reasons similar to those mentioned in the above sections. It should be noted that in this case, just of the same way as above, there appears a time $T_3 > 0$, until this time, the smoothness of the solution will continue.

## 6. Another case of smoothness at variable $t$.

Assume the dates of the studied Problem 1(1.1)-(1.3) are sufficiently smooth. If this problem has a solution, then it must also possess such smoothness, and we can differentiate both parts of the equation from this problem. Then, differentiating equation (1.1) on $t$ we get to the following equation

(6.1) $$u_{tt} - \nu\Delta u_t + \big((u \cdot \nabla)u\big)_t + \nabla p_t = f_t$$

So we will study the problem (6.1), (1.2), and with the initial conditions $u(0) = u_0$, $u_t(0) = u_1$.

In addition, we will assume the solution is bounded, i.e., there exists a constant $K > 0$ such that $\|u\| \leq K$, which follows from the previous results. We will act in the following way: here we will use the following operator

$$A_1(D_t, D_x)u = u_{tt} - \Delta u_t + u_t + u$$

Thus, acting as in the previous cases, we obtain the following equation

$$0 = \langle u_{tt} - \nu\Delta u_t + \nabla p_t, u_{tt} - \Delta u_t + u_t + u \rangle +$$
$$+ \langle \big((u \cdot \nabla)u\big)_t - f, u_{tt} - \Delta u_t + u_t + u \rangle =$$

(6.2) $$= \|u_{tt}\| + \frac{1}{2}\frac{d}{dt}\big(\int |\nabla u_t|^2 + \nu|\Delta u_t|^2 + \nu|\nabla u_t|^2\big)dx +$$
$$+ \nu\|\nabla u_t\|^2 + \langle u_{tt}, u \rangle - \int f \times (A_1(D_t, D_x)u)dx +$$
$$+ \int \{[D_t\big((u \cdot \nabla)u\big) + \nabla p_t] \times (A_1(D_t, D_x)u)\}dx$$

Whence,

$$\langle \nabla p_t, A_1(D_t, D_x)u \rangle = \langle p, (D_t^2 - D_t\Delta + D_t + I)\nabla u = 0 \rangle,$$

due to the conditions of the studied Problem.

In the equality (6.2), one must investigate firstly the following term

$$\int \big[D_t\big((u \cdot \nabla)u\big) \times A_1(D_t, D_x)u\big] dx.$$

If we estimate each term from the above expression, then we get

$$\int \big[|D_t((u \cdot \nabla)u)u_{tt}| + |D_t((u \cdot \nabla)u)\Delta u_t|\big]dx \leq \frac{4\nu + 2}{\nu}\int \big(D_t((u \cdot \nabla)u)\big)^2$$
$$+ \frac{1}{4}\int (u_{tt})^2 dx + \frac{\nu}{2}\int (\Delta u_t)^2 dx,$$

$$\int \big|\big(D_t((u \cdot \nabla)u)\big)u_t\big| dx$$
$$+ \int \big|\big(D_t((u \cdot \nabla)u)\big)u\big| dx$$
$$\leq \frac{1}{4}\int \big(D_t((u \cdot \nabla)u)\big)^2 dx + \frac{1}{4}\int ((u_t)^2 + u^2)dx$$

Thus, we obtain the following inequality

$$\int |D_t((u \cdot \nabla)u) \times (u_{tt} - \Delta u_t + u_t + u)|dx \leq \frac{4\nu + c}{\nu}\left(\int \big(D_t((u \cdot \nabla)u)\big)^2\right)dx$$
$$+ \frac{1}{8}\big(\|u_{tt}\|^2 + \nu\|\Delta u_t\|^2 + \|u_t\|^2 + \|u\|^2\big) \quad (6.3)$$

It is not difficult to see that the following terms on the right side of the inequality mentioned above

$$\frac{1}{4}\|u_{tt}\|^2 + \frac{\nu}{2}\|\Delta u_t\|^2 + \frac{1}{4}\|u_t\|^2 + \frac{1}{8}\|u\|^2$$

we can estimate by using norms of terms from inequality (6.2).

So, it is necessary to estimate the expression $\int \left(D_t((u \cdot \nabla)u)\right)^2 dx$, which we can act in the following way

$$\int \left(D_t((u \cdot \nabla)u)\right)^2 dx = \int \left(((u_t \cdot \nabla)u) + ((u \cdot \nabla)u_t)\right)^2 dx$$
$$\leq \frac{\nu + c}{4\nu} \int (|u_t|^2 + \nu|\nabla u|^2 + |u|^2 + (1+\nu)|\nabla u_t|^2)^2 dx \quad (6.4)$$

Taking into account the obtained inequalities (6.3) and (6.4) in equation (6.2), and carrying on some calculations, we arrive

$$\frac{1}{2}\|u_{tt}\|^2 + \frac{\nu}{2}\|\Delta u_t\|^2 + \nu\|\nabla u_t\|^2 + \frac{1}{2}\left|\int D_t((1+\nu)|\nabla u_t|^2 + \nu|\nabla u|^2 + |u_t|^2 + |u|^2)\right| dx$$
$$\leq \frac{4(\nu + c_1)}{\nu} \int (|\nabla u_t|^2 + \nu|\nabla u|^2 + |u_t|^2 + |u|^2)^2 dx$$
$$+ \frac{(c_2 + \nu)}{\nu}\|f\|^2_{W^{2,2}(H) \cap W^{1,2}(H^2)} + CM(\cdots) + c_3\|u_0\|\|u_1\|$$

where $c_1, c_2, c_3$ and $C$ are constants independent of $u$; here $K$ is the constant defined at the beginning of this section.[4]

As above, we arrive at the result that we must study the appropriate Cauchy problem. So, it needs to study the following Cauchy problem for ordinary differential inequation

$$\frac{d}{dt}(|\nabla u_t|^2 + \nu|\nabla u|^2 + |u_t|^2 + |u|^2)$$
$$\leq C_0(\nu, \varepsilon, M)(|\nabla u_t|^2 + \nu|\nabla u|^2 + |u_t|^2 + |u|^2)$$
$$+ C_1(\nu, \varepsilon, M, \|f\|) \quad (6.5)$$

with the initial conditions, which follow from the initial conditions posed at the beginning of this Section

$$u(0,x) = u_0(x), u_t(0,x) = u_1(x), \nabla u(0,x) = \nabla u_0(x), \nabla u_t(0,x) = u_1(x)$$

according to the inequalities (6.3) and (6.4).

As in the previous sections, the obtained Cauchy problem will allow us to receive the necessary a priori estimates for the considered problem. So, by carrying out analogous calculations as in the previous sections, one can receive inequalities similar to inequalities of such types from the above sections.

Consequently, by repeating the arguments and by carrying out the corresponding calculations similar to those leading to the proof of the results mentioned above, we obtain the necessary estimations. For brevity, we will not conduct these calculations; we only wish to note that, as in all the above sections, also here, there exists a time $T_4 > 0$ of the blowup that bounds the time $T < T_4$ of the smoothness of the solutions.

The results obtained in this article demonstrate that by continuing to apply this method, one can show other properties related to the smoothness of the solution to the Navier-Stokes equation in the 3-dimensional case.

---

[4] We should note that if $u \in W^{2,2}(R_+; H^2)$ then

$$\|u_t\|^2_{L^2(H)} \leq \frac{\varepsilon}{2}\|u_{tt}\|^2 + \frac{1}{2\varepsilon}\|u\|^2_{L^2(H)} + u_t(0,x)u(0,x).$$

References

[1] Caffarelli L., Kohn R., Nirenberg L., Partial regularity of suitable weak solutions of the Navier-Stokes equation, Comm. Pure Appl. Math., 35 (1982), 771-831

[2] Caffarelli L., A. Vasseur A., The De Giorgi method for nonlocal fluid dynamics, Nonlinear PDE Advanced Courses in Mathematics, CRM, Barcelona (2012), 1-38

[3] Chamorro D., Lemarié-Rieusset P-G., Mayoufi K., The Role of the Pressure in the Partial Regularity Theory for Weak Solutions of the Navier–Stokes Equations, Arch. Rational Mech. Anal. 228, 237–277 (2018).

[4] Constantin P., Euler Equations, Navier–Stokes Equations and Turbulence, Departments of Mathematics, C.I.M.E.Lectures, 2/21/2004, Princeton.

[5] Diestel, J., Geometry of Banach spaces, Lecture Notes in Math, 485, Springer-Verlag, 1975.

[6] Dong H., D. Du D., Partial regularity of weak solutions of the Navier–Stokes equations in at the first blow-up time, Comm. Math. Phys., 273 (3) (2007), 785-801.

[7] Fefferman Ch. L., Existence and Smoothness of the Navier–Stokes Equation, Millennium Problems, Clay Math. Inst., (2006).

[8] Galdi G., An Introduction to the Mathematical Theory of the Navier–Stokes Equations: Steady-State Problems, (second edition), Springer (2011)

[9] Hopf E., ̈Uber die anfangswertaufgabe für die hydrodynamischen grundgleichungen, Math. Nachr., 4, 213-231(1951)

[10] Hopf E., On nonlinear partial differential equations, Lecture Series of the Symposium on Partial Differential Equations, Berkeley, 1955, Ed. Univ. Kansas, I-29, (1957).

[11] Kukavica I., On partial regularity for the Navier–Stokes equations. Discrete Contin.Dyn. Syst. 21, 717–728, (2008).

[12] Ladyzenskaja O. A., Seregin G., On partial regularity of suitable weak solutions to the threedimensional Navier–Stokes equations, J. Math. Fluid Mech., 1 (4) (1999), 356-387 [13] Leray J., Sur le Mouvement d'un Liquide Visquex Emplissent l'Espace, Acta Mat h. J. 63, 193–248, (1934).

[14] Lin F., A new proof of the Caffarelli–Kohn–Nirenberg Theorem, Comm. Pure Appl. Math., 51 (1998), 241-257

[15] Lindberg S., On the integrability properties of Leray-Hopf solutions of the Navier-Stokes equations on $R^3$, arXiv preprint arXiv:2412.13066, 2024

[16] Lions J.-L., Quelques methodes de resolution des problemes aux limites nonlineares, Dunod, Gauthier-Villars, Paris, (1969).

[17] Scheffer V., Partial regularity of solutions to the Navier–Stokes equations, Pacific J. Math., 66 (2) (1976), 535-552

[18] Serrin J., On the interior regularity of weak solutions of the Navier–Stokes equations Arch. Ration. Mech. Anal., 9 (1962), 187-195

[19] Seregin, G.A., Differentiability properties of weak solutions to te Navier-Stokes equations, Algebra and Analysis, 14(2002), No. 1, 193-237

[20] Soltanov K. N., Remarks on the uniqueness the weak solutions to the incompressible Navier– Stokes-equations, Russian Journal of Mathematical Physics, 31, 3,(2024)

[21] Soltanov K. N., *Some applications of the nonlinear analysis to the differential equations*, Baku, ELM, (in Russ.), (2002).

[22] Soltanov K. N., Nonlinear Operators, Fixed-Point Theorems, Nonlinear Equations, "Current Trends in Analysis and Its Applications", Conference on Trends in Mathematics 2013Krakow,(2015) DOI 10.1007/978-3-319-12577-0 41

[23] Soltanov K. N., On semi-continuous mappings, equations and inclusions in a Banach space, Hacettepe Journal of Math and Statistics, V. 37 (1) (2008)

[24] Soltanov K. N., Perturbation of the mapping and solvability theorems in Banach space, Nonlinear Analysis, (2010), doi:10.1016/j.na.2009.06.067

[25] Soltanov K. N. , Sprekels J., Nonlinear equations in nonreflexive Banach spaces and fully nonlinear equations. Advances in Mathematical Sciences and Applications, 1999, v.9, n 2, 939-972.

[26] Soltanov K. N., Ahmadov M., Solvability of Equation of Prandtl-von Mises type, Theorems of Embedding. Trans. NAS Azerb, ser. Phys.-Tech. & Math., 37, 1, (2017).

[27] Soltanov K. N., Some embedding theorems and nonlinear differential equations. Trans. Ac. Sci. Azerb., ser. Math. & Phys.-Tech., 19, 5, 125 - 146, (1999).

[28] Struwe M., Partial regularity results for the Navier–Stokes equations, Comm. Pure Appl. Math., 41 (1988), 437-458

[29] Tao T., Finite time blowup for an averaged three-dimensional Navier-Stokes equation, arXiv:1402.0290v3 [math.AP] 1 Apr 2015.

[30] Temam R., Navier–Stokes Equations Theory and Numerical Analysis, North-Holland Pub. Comp., in English, 3rd rev. edn, (1984).

[31] Vasseur A., A new proof of partial regularity of solutions to Navier–Stokes equations, Nonlinear Differential Equations Appl., 14 (2007), 753-785

[32] Wojciech S. Oż˙an´ski, Partial regularity of Leray–Hopf weak solutions to the incompressible Navier–Stokes equations with hyperdissipation, Analysis and Partial Differential Equations 16(3):747-783, (2023)

[33] Yanqing Wang, Gang Wu, A unified proof on the partial regularity for suitable weak solutions of non-stationary and stationary Navier–Stokes equations, Journal of Differential Equations, 256, 3, (2014), 1224-1249

National Academy of Sciences of Azerbaijan, Baku, AZERBAIJAN
*Email address*: sulta.kamal.n@gmail.com
*URL*: https://www.researchgate.net/profile/Kamal-Soltanov/research